\documentclass[letter,12pt]{article}
\usepackage{amsfonts}
\usepackage{amssymb}
\usepackage{amsmath}
\usepackage{cite}
\newtheorem{theo}{Theorem}

\newtheorem{prop}[theo]{Proposition}
\newtheorem{lema}[theo]{Lemma}

\newtheorem{obs}[theo]{Remark}

\newcommand{\R}{\mathbb{R}}

\newcommand{\C}{\mathbb{C}}

\newcommand{\Ti}{\mathcal{T}}

\title{On the projective derivative cocycle\\ 
for circle diffeomorphisms}
\author{Andr\'es Navas \, and \, Mario Ponce}
\begin{document}
\maketitle

\vspace{-0.5cm}

\noindent{\bf Abstract.} We study the projective derivative as a cocycle of M\"obius transformations over groups of circle diffeomorphisms. 
By computing precise expressions for this cocycle, we obtain several results about reducibility and almost reducibility to a cocycle of 
rotations. We also introduce an extension of this cocycle to the diagonal action on the 3-torus for which we generalize the previous results. 

\vspace{0.2cm}


\vspace{0.2cm}


\vspace{0.9cm}


\noindent {\Large{\bf Introduction}}
\vspace{0.5cm}

A classical and useful idea to study a given dynamical system consists in looking for approximations that are easier to handle but still reveal 
some of its relevant features. For differentiable systems, it is worth to look at the first-order reduction, that is, the classical derivative. 
This naturally becomes a cocycle over the action, the cocycle relation being just the chain rule. A fundamental object  to deal with in this 
context are the Lyapunov exponents, and a relevant feature to study is the (almost) reducibility of the cocycle.

By definition, cocycles take values in groups, the group associated to the derivative cocycle just being the linear group. Now, the existence of larger 
but still finite dimensional groups of transformations naturally give raise to other cocycles, which induce approximations of systems in higher regularity. 
In this note we deal with a remarkable cocycle of this type, which arises from the projective derivative. Given a $C^2$ circle diffeomorphism, this corresponds 
to the map that sends each point of the circle to the unique M\"obius transfiormation whose 2-jet at this point coincides with that of the diffeomorphism. 
This idea can be traced back to \cite{SING}, although part of the attribution should be deserved to Thurston, 
according to certain sources \cite{GALA}. 

The projective derivative can be seen as a cocycle taking values in the group $\mathcal{M}$ of M\"obius transformations, which we identify 
with $\mathrm{PSL}(2,\mathbb{R})$. This cocycle can be considered over single circle diffeomorphisms and, more generally, over groups of 
such diffeomorphisms. We establish several basic properties and two main results. The first of these (Theorem A below) is a characterization 
of those groups for which the projective cocycle is $\mathrm{SO}(2,\mathbb{R})$-reducible, that is, conjugate (cohomologous) to a cocycle of 
rotations. The second one (Theorem B) is a characterization of groups for which the cocycle is almost $\mathrm{SO}(2,\mathbb{R})$-reducible. 
In the appendix, we extend the definition of this cocycle to the 3-torus (for the diagonal action), and we provide generalizations of our results for 
the projective cocycle to this broader context. All the terminology is carefully explained throughout the text. 


\section{The projective derivative of a $C^2$ circle diffeomorphism}

Troughout this work, all maps will be assumed to be orientation preserving.

Following \cite{SING}, given a $C^2$ circle diffeomorphism $f$, the {\em projective derivative} 
of $f$ at a point $\theta$ is the unique M\"obius transformation $P = P_{f,\theta}$ of the circle that coincides up to 
order 2 with $f$ at $\theta$, that is, 
$$P (\theta) = f(\theta), \qquad D P (\theta) = Df (\theta), \qquad D^2 P (\theta) = D^2 f (\theta).$$
Below we proceed to explicitly compute this map. Since we use complex notation\footnote{In some places this 
simplifies computations, in some others this makes it more involved. The reader is invited to carry out his/her 
own computations in the coordinates he/she prefers (computations related to the Appendix are particularly 
enlightening in projective coordinates).}, we identify an 
angle $\theta \in [0,2\pi]$ with the point $e^{i \theta} \in \mathbb{S}^1$. Looking at the action of 
$P = P_{f,\theta}$ on the boundary circle $\mathbb{S}^1$, we have 
\begin{equation}\label{e:gen-normal}
e^{i P(\xi)} = e^{i\kappa} \cdot \frac{e^{i \xi} - \sigma}{1 - \bar{\sigma} e^{i\xi}},
\end{equation}
where $\kappa = \kappa(f,\theta)$ lies in $[0,2\pi]$ and $\sigma = \sigma(f,\theta)$ 
lies in the Poincar\'e disk $\mathbb{D}$. This yields 
$$P(\xi) = 
\frac{1}{i} \log \left[  e^{i\kappa} \cdot \frac{e^{i \xi} - \sigma}{1 - \bar{\sigma} e^{i\xi}} \right] 
= \kappa + \frac{1}{i} \log (e^{i \xi} - \sigma) - \frac{1}{i} \log(1 - \bar{\sigma} e^{i \xi}).$$
Hence,
$$D P  (\xi) 
=
\frac{1}{i} \left[ \frac{i e^{i\xi}}{e^{i\xi} - \sigma} 
+ \frac{ie^{i\xi} \bar{\sigma}}{1 - \bar{\sigma}e^{i\xi}} \right]
= \frac{1}{i} \left[ \frac{ie^{i\xi}(1 - \bar{\sigma}e^{i\xi} + \bar{\sigma}e^{i\xi} - \sigma \bar{\sigma})}{(e^{i\xi}-\sigma)(1-\bar{\sigma}e^{i\xi})} \right]
= \frac{1 - \|\sigma\|^2}{\|1 - \bar{\sigma}e^{i\xi}\|^2}.$$
Moreover,
\begin{eqnarray*}
D^2 P (\xi)  
&=& 
(1 - \|\sigma\|^2) \frac{d}{d\xi} \left[ \frac{1}{(1 - e^{-i\xi}\sigma)(1 - e^{i\xi}\bar{\sigma})} \right]\\
&=&
(1-\|\sigma\|^2)(-1) \left[ \frac{ie^{-i\xi}\sigma(1-e^{i\xi}\bar{\sigma}) - 
i \bar{\sigma} e^{i\xi} (1 - {\sigma} e^{-i\xi})}{(1 - e^{-i\xi} \sigma)^2 
(1 - e^{i\xi} \bar{\sigma})^2} \right]\\
&=&
\frac{(1-\|\sigma\|^2) (-1) i (e^{-i\xi} \sigma  - \overline{e^{-i\xi} \sigma})}{\|1-\bar{\sigma}e^{i\xi}\|^4}\\
&=&
\frac{(1-\|\sigma\|^2) \, 2 \, \mathrm{Im} (e^{-i\xi}{\sigma})}{\|1-\bar{\sigma}e^{i\xi}\|^4}.
\end{eqnarray*}

Summarizing:
$$P(\xi) = 
\frac{1}{i} \log \left[  e^{i\kappa} \cdot \frac{e^{i \xi} - \sigma}{1 - \bar{\sigma} e^{i\xi}} \right],$$
$$D P (\xi) = \frac{1 - \|\sigma\|^2}{\|1 - \bar{\sigma}e^{i\xi}\|^2},$$
$$D^2 P (\xi) = \frac{(1-\|\sigma\|^2) \, 2 \, \mathrm{Im} (e^{-i\xi}{\sigma})}{\|1-\bar{\sigma}e^{i\xi}\|^4}.$$
Let us denote $\mathrm{D} := Df (\theta)$ and $\Delta := D^2 f (\theta)$. To simplify notation, we also let 
$$\bar{\sigma} e^{i\theta} = r e^{i \beta}, \qquad 0 \leq r < 1.$$
Then, the first order coincidence between $f$ and $P$ at $\xi=\theta$ becomes
$$\mathrm{D} = \frac{1-r^2}{1 - 2 r \cos(\beta) + r^2},$$
hence
$$\cos (\beta) = \frac{\mathrm{D} (1+r^2)- (1-r^2)}{2 \, \mathrm{D}\, r}.$$
The second order relation  becomes
$$\Delta = -\frac{(1 - r^2) \, 2 r \sin(\beta)}{[1 - 2r \cos(\beta) + r^2]^2},$$
thus
$$\sin (\beta) =  -\frac{\Delta}{\mathrm{D}^2} \left( \frac{1 - r^2}{2r} \right).$$
Using that \, $\sin^2 (\beta) + \cos^2 (\beta) = 1$, \, we obtain
$$\left[ \mathrm{D}(1+r^2) - (1-r^2) \right]^2 + \left( \frac{\Delta \, (1-r^2)}{\mathrm{D}} \right)^2 
= 4 \, \mathrm{D}^2 \, r^2,$$
which is equivalent to the equation
$$r^4 (\mathrm{D}^2(\mathrm{D}+1)^2 + \Delta^2) - 2r^2 (\mathrm{D}^2(\mathrm{D}^2+1)+\Delta^2) 
+ \mathrm{D}^2(\mathrm{D}-1)^2 + \Delta^2 = 0.$$
This equation has two solutions, namely $r^2 = 1$ and 
\begin{equation}\label{r_sub_ene}
r^2 = \frac{\mathrm{D}^2(\mathrm{D}-1)^2 + \Delta^2}{\mathrm{D}^2(\mathrm{D}+1)^2 + \Delta^2}.
\end{equation}
The former is impossible as $\sigma$ lies inside the unit disc 
$\mathbb{D}$ and $\| \sigma \| = r$, so that we necessarily have  
\begin{equation}\label{erre}
r = \sqrt{ \frac{\mathrm{D}^2(\mathrm{D}-1)^2 + \Delta^2}{\mathrm{D}^2(\mathrm{D}+1)^2 + \Delta^2}}.
\end{equation}
Consequently,
$$\sin (\beta) = \frac{-2\mathrm{D}\Delta}{\sqrt{ (\mathrm{D}^2(\mathrm{D}-1)^2 + \Delta^2) (\mathrm{D}^2(\mathrm{D}+1)^2 + \Delta^2) }},$$
$$\cos (\beta) = \frac{\mathrm{D}^2 (\mathrm{D}^2 - 1) + \Delta^2}{\sqrt{ (\mathrm{D}^2(\mathrm{D}-1)^2 + \Delta^2) (\mathrm{D}^2(\mathrm{D}+1)^2 + \Delta^2) }}.$$   
We can also recover the value of $\sigma=re^{i\alpha}$ using 
\begin{equation}\label{e:angle-sigma}
\alpha = \theta-\beta, \qquad  \tan (\beta) = \frac{-2\mathrm{D}\Delta}{\mathrm{D}^2 (\mathrm{D}^2 - 1) + \Delta^2}.
\end{equation}
Notice that these values are not well defined when $\mathrm{D}=1$ and $\Delta = 0$. However, in this case, (\ref{erre}) implies that $\sigma = 0$. 
Finally, the value of $\kappa$ can be recovered from the zero order relation between $f$ and $P$ at $\theta$:
$$f(\theta)=\kappa+\frac{1}{i}\log \left[\frac{e^{i\theta}-\sigma}{1-\overline{\sigma}e^{i\theta}}\right].
$$


\section{The projective cocycle and its reduction} 

It readily follows from the definition that the projective derivative is a cocycle, which means that for all $C^2$ circle diffeomorphisms 
$f,g$ and every $\theta$, one has
$$P_{f,g(\theta)} \circ P_{g , \theta} = P_{fg, \theta}.$$
Consequently, one also has
$$P_{f^{-1},\theta} = P_{f, f^{-1}(\theta)}^{-1}.$$
In particular, given an action of a group by $C^2$ circle diffomorphisms, the map $f \to P_{f,\cdot}$ is a cocycle above this action taking 
values in the M\"obius group $\mathcal{M}$.

Recall that, given a cocycle $P$ over a group action on a space $X$ and taking 
values in a group $G$, one says that $P$ is $H$-reducible for a subgroup 
$H$ of $G$ if it is cohomologous to a cocycle taking values in $H$. This means that  
there is a continuous map $B: X \to G$ such that $B(f(x))\cdot P (x)\cdot B(x)^{-1}$ belongs to $H$ for every $x\in X$ and each $f$ in the 
acting group $\Gamma$. Below we are interested in the case where $P$ is the projective cocycle (hence $G = \mathcal{M}$) and $H$ is the group of rotations.

\vspace{0.5cm}

\noindent{\bf Theorem A} 
{\em Given a group of $C^2$ circle diffeomorphisms, the projective cocycle associated to it 
is $\mathrm{SO} (2,\mathbb{R})$-reducible if and only if the action is $C^2$ conjugate to an action by rotations.}

\vspace{0.5cm}

The direct implication of this theorem is immediate. Indeed,   
assume that the group action is $C^2$ conjugate to an action by rotations, and denote by 
$\varphi$ any such conjugating map. This means that, for all $f$ in the acting group $\Gamma$, we have 
$$\varphi \circ f \circ \varphi^{-1} = R_{\rho(f)},$$
where $\rho(f)$ denotes the rotation number of $f$. Since \, $P_{R,\cdot} = R$ \, obviously holds for every rotation 
$R$ (and, in general, for every M\"obius transformation), using the cocycle relation of $P$ in the last equality we obtain 
$$P_{\varphi, f \varphi^{-1}(\theta)} \circ P_{f, \varphi^{-1} (\theta)} \circ P_{\varphi^{-1} , \theta} = P_{R_{\rho(f)}, \theta} = R_{\rho(f)}.$$
Using now the relation $P_{\varphi^{-1}, \theta} = P_{\varphi, \varphi^{-1} (\theta)}^{-1}$ and then changing the $\theta$ by 
$\varphi (\theta)$, this reduces to 
\begin{equation}\label{eq:conjj}
P_{\varphi, f(\theta)} \circ P_{f,\theta} \circ P_{\varphi,\theta}^{-1} = R_{\rho(f)}.
\end{equation}
Therefore, the map $B$ defined by $B(\theta) := P_{\varphi, \theta}$ performs the desired reduction of $P$. 

\vspace{0.2cm}

The reverse implication is much subtler. Assume that $P$ is cohomologous to a cocycle with values in $\mathrm{SO}(2,\mathbb{R})$ 
via a (continuous) map $B$. We first claim that, in this case, the derivatives and second derivatives of group elements are uniformly 
bounded. Indeed, from
\begin{equation}\label{reduction}
B (f(\theta)) \circ P_{f , \theta} \circ B (\theta)^{-1} = R_{f,\theta}
\end{equation}
we obtain 
$$P_{f, \theta} = B (f(\theta))^{-1} \circ R_{f,\theta} \circ B (\theta).$$
Since $B$ takes values in a compact subset of $\mathcal{M}$ and $R_{f,\theta}$ belongs to $\mathrm{SO}(2,\mathbb{R})$, we conclude that 
$P$ takes values in a compact subset of $\mathcal{M}$. Now, elements of $\mathcal{M}$ can be expressed in the form
$$z \mapsto e^{i\kappa} \cdot \frac{z - \sigma}{1-\bar{\sigma} z}, \qquad \kappa \in [0,2\pi], \quad \sigma \in \mathbb{D}.$$
A subset of $\mathcal{M}$ is relatively compact if and only if it is contained in a set of maps as above for which there is a uniform 
bound $\lambda < 1$ for the values of $\|\sigma\|$ arising along this set. In particular, if $\lambda$ is such an upper bound 
for the projective cocycle $P$, relation (\ref{erre}) yields, for all $f \in \Gamma$ and all $\theta$, 
$$\sqrt{\frac{\mathrm{D}^2(\mathrm{D}-1)^2 + \Delta^2}{\mathrm{D}^2(\mathrm{D}+1)^2 + \Delta^2}} \leq \lambda,$$
where $\mathrm{D} = Df (\theta)$ and $\Delta = D^2 f (\theta)$. Elementary computations then show that 
$$\mathrm{D} \leq \frac{1+\lambda}{1-\lambda}, \qquad \Delta \leq \frac{2 \, (1 + \lambda)}{(1 - \lambda)^2},$$
thus closing the proof of uniform boundedness of derivatives and second derivatives of group elements.

\vspace{0.2cm}

To pursue the proof of the conjugacy, we next state and prove a proposition of independent interest.

\vspace{0.2cm} 

\begin{prop}\label{independent} Let $\Gamma$ be a group of $C^1$ circle diffeomorphisms. 
If there is a uniform bound for the derivatives of all group elements, then $\Gamma$ is 
conjugate to a group of rotations by a $C^1$ diffeomorphism.
\end{prop}

\vspace{0.2cm}

This corresponds to a version of Herman's criterium in \cite[Chapitre IV]{herman}  (see also \cite[\S 3.6.2]{book}).
However, he only deals with single diffeomorphisms with irrational rotation number, and some care is needed for 
general group actions. We just sketch the complete argument since, although the claim is not explicitly stated 
in the literature, it is well known to the specialists. 

\vspace{0.5cm}

\noindent{\bf Proof of Proposition \ref{independent}.} 
Since derivatives are uniformly bounded, passing to the inverse group elements we deduce that they are also bounded 
away from zero. This easily implies that the action is equicontinuous and, using Ascoli-Arzela's theorem, this yields that 
it is topologically conjugate to an action by rotations \cite[\S 2.3]{book}. As a consequence, $\Gamma$ is Abelian. 

We claim, however, that the conjugacy map can be taken to be bi-Lipschitz. This can be established via the cohomological equation 
associated to the conjugacy problem. However, a more concrete argument (which will be useful later) proceeds as follows. 
Assume for a while that $\Gamma$ is finitely generated. Following \cite{compositio}, we let 
\begin{equation}\label{conjugadores}
\tilde{\varphi}_n (\theta) = \frac{1}{B(n)} \sum_{f \in B(n)} \tilde{f} (\theta),
\end{equation}
where $B(n)$ denotes the ball of radius $n$ in the group, and the tilde stands for appropriate lifts of maps. 
This induces a circle diffeomorphism $\varphi_n$ for which the sequence of conjugates 
$$\varphi_n \circ f \circ \varphi_n^{-1}$$
converge to $R_{\rho(f)}$ for every group element $f$. It follows from the hypothesis and the definition that there is a uniform 
bound for the derivatives of the maps $\varphi_n$, and that these are also bounded away from zero. Therefore, up to a subsequence, 
they converge to a bi-Lipschitz circle homeomorphism that conjugates $\Gamma$ to a group of rotations. In case the group 
is not finitely generated, one can apply this argument to every finitely generated subgroup. Since the upper and lower bounds 
for the derivatives of the conjugating maps are independent of the subgroup, a standard compactness argument closes 
the proof of the existence of a bi-Lipschitz conjugating map.

Assume next that the group is finite. Since it is isomorphic to a group of rotations, it is cyclic, say generated by a diffeomorphism 
$f$ of order $k$. Let $\tilde{f}$ be a lift of $f$ satisfying $\tilde{f}^k (\tilde{\theta}) = \tilde{\theta} + 2\pi k \varrho(\tilde{f})$ 
for all $\tilde{\theta}$, where $\varrho(\tilde{f})$ is the translation number of $\tilde{f}$. The map 
\begin{equation}\label{caso-finito}
\tilde{\varphi} (\tilde{\theta}) : = \frac{\tilde{\theta} + \tilde{f} (\tilde{\theta}) + \ldots + \tilde{f}^{k-1}(\tilde{\theta})}{2 \pi k}
\end{equation}
is a $C^1$ diffeomorphism of the real line that satisfies \, $\tilde{\varphi} ( \tilde{f}(\tilde{\theta})) = \tilde{\varphi}(\tilde{\theta}) + \varrho (\tilde{f})$ \, 
for all $\tilde{\theta}$. Besides, $\tilde{\varphi}$ commutes with the translation $\tilde{\theta} \to \tilde{\theta}+2\pi$, hence induces a $C^1$ circle 
diffeomorphism $\varphi$. Finally, the previous relation descends into 
\, $\varphi (f( \theta )) = R_{\rho(f)} (\varphi(\theta)),$ \, which is the desired conjugacy relation.

Assume now that the group is infinite. Since it is bi-Lipschitz conjugate to a dense group of rotations, its action on the circle 
is ergodic with respect to the Lebesgue measure. A direct application of \cite{3-remarks} (see also \cite[\S 3.6.2]{book})  
then yields that every bi-Lipschitz conjugacy to a group of rotations is, actually, a $C^1$ diffeomorphism. 
$\hfill\square$

\vspace{0.4cm}

Recall that we are dealing with a group for which not only the first derivatives but also the second ones are uniformly 
bounded. In this framework, a simple argument proves the next lemma.





\vspace{0.2cm}

\begin{lema} If $\, \Gamma$ is a group of $C^2$ circle diffeomorphisms for which derivatives and second derivatives are 
uniformly bounded along group elements, then $\Gamma$ is conjugate to a group of rotations by a $C^{1+Lip}$ diffeomorphism.
\end{lema}

\noindent{\bf Proof.} If $\Gamma$ is finitely generated, we can use the sequence of conjugating diffeomorphisms (\ref{conjugadores}).  
(In case of a finite group, we just use (\ref{caso-finito}).) It readily follows from the definition that these 
have uniformly bounded first and second derivatives. By the Ascoli-Arzela's theorem, any limit of them will be of class $C^{1+Lip}$. 
In case of non finitely generated groups, a compactness argument as above gives the desired $C^{1+Lip}$ conjugating diffeomorphism. 
$\hfill\square$

\vspace{0.4cm}

Passing from a $C^{1+Lip}$ to a $C^2$ conjugacy is quite tricky. In \cite[Chapitre IV]{herman}, Herman proposes a proof that works 
for single diffeomorphisms with irrational rotation number; this  can be pursued general groups using the arguments of \cite{3-remarks}. 
However, for the proof of Theorem~A, we prefer to give a simpler argument that uses some of the information contained in the 
projective cocycle. Let hence $\varphi$ be a $C^{1+Lip}$ circle diffeomorphism that conjugates our group to a group of rotations. Taking 
the affine derivative on both sides of the relation
$$\varphi \circ f \circ \varphi^{-1} = R_{\rho(f)},$$
we obtain that the following equality holds for almost every $\theta$:
\begin{equation}\label{eq:u}
\frac{D^2 \varphi}{D \varphi} (f(\theta)) \cdot Df (\theta) + \frac{D^2 f}{D f} (\theta) - \frac{D^2 \varphi}{D \varphi} = 0.
\end{equation}
Notice that the function $u := D^2 \varphi / D \varphi$ is known only to belong to $L^{\infty}$. To conclude the proof of 
Theorem A, we need to prove that it 
coincides with a continuous function. To do this, we come back to the equation of reduction of the projective cocycle. Namely, taking 
the affine derivative at both sides of the equation (\ref{reduction}) we obtain, for each $\theta$ and every point $\xi$, 
$$\frac{D^2 B (f(\theta))}{D B (f(\theta))} (P_{f,\theta} (\xi)) \cdot DP_{f,\theta} (\xi) 
+ \frac{D^2 P_{f,\theta} (\xi)}{D P_{f,\theta} (\xi)} - \frac{D^2 B(\theta)}{D B(\theta)} (\xi) = 0,$$
where derivatives are taken with respect to $\xi$. Evaluating this equation at $\xi = \theta$ and using the properties of $P$ then yields
$$\frac{D^2 B (f(\theta))}{D B (f(\theta))} ( f (\theta)) \cdot Df (\theta) 
+ \frac{D^2 f}{D f} (\theta) - \frac{D^2 B(\theta)}{D B(\theta)} (\theta) = 0.$$
Let us consider the function 
$$v(\theta) := \frac{D^2 B(\theta)}{D B(\theta)} (\theta).$$
This is a continuous function, and the previous relation together with (\ref{eq:u}) read as 
$$\frac{D^2 f }{Df} = v - (v \circ f) \cdot D f= u - (u \circ f) \cdot D f.$$
As a consequence, 
$$(u-v) \circ f \cdot Df = u-v.$$
This easily implies that the function $u-v$ is the density of a (signed, but finite) measure that is invariant under the group action. 
Now, if the group is infinite, then it has a unique invariant probability measure, which is the push-forward of the Lebesgue 
measure by the conjugating map. Since this map is $C^1$, this probability measure has a continuous density. We hence 
conclude that $u-v$, which is a scalar multiple of this density, is continuous. As a consequence, $u$ is continuous, thus 
concluding the proof of Theorem A  for infinite groups. For finite groups, just follow the argument using (\ref{caso-finito}).

\vspace{0.2cm}

\begin{obs} \label{rem:bounded}
{\em Except for the last argument, in the proof of Theorem A we just used the fact that $P$ is a uniformly bounded cocycle (as a function defined on 
the acting group and the circle), where $P$ is thought as a matrix cocycle via the identification of $\mathcal{M}$ with $\mathrm{PSL}(2,\mathbb{R})$. 
We claim, however, that this condition is equivalent to $\mathrm{SO}(2,\mathbb{R})$-reducibility. Indeed, by \cite{CNP}, boundedness of $2 \times 2$ 
matrix cocycles and $\mathrm{SO}(2,\mathbb{R})$-reducibility are equivalent in the case of a minimal dynamics on the basis. Now, for the projective 
cocycle, boundedness also implies that the group is topologically conjugate to a group of rotations (see the proof of Proposition 1, where actually 
$C^1$ conjugacy is established). Therefore, the minimality hypothesis is guaranteed for infinite groups; in this case, the result of \cite{CNP} 
yields $\mathrm{SO}(2,\mathbb{R})$-reducibility, and Theorem A applies. For finite groups, again just follow the argument using (\ref{caso-finito}).}
\end{obs}


\section{Almost reduction of the projective cocycle} 

Recall that, given a cocycle $P$ over a group action on a space $X$ taking values in a group $G$ (with a prescribed 
left-invariant metric), one says that 
$P$ is almost $H$-reducible if there exists a sequence of cohomologous cocycles converging to maps into $H$. More 
precisely, there exists a sequence of continuous maps $B_n: X  \to G$ such that, for all $x \in X$ and all 
$f$ in the acting group $\Gamma$, the distance between 
$$B_n (f(x)) \cdot P_{f,x} \cdot B_n (x)^{-1}$$ 
and $H$ uniformly converges to zero.

\vspace{0.5cm}

\noindent{\bf Theorem B} 
{\em If $\Gamma$ is a finitely generated Abelian group of $C^2$ circle diffeomorphisms acting 
freely, then the projective cocycle above its action is almost $\mathrm{SO}(2,\mathbb{R})$-reducible.}

\vspace{0.5cm}

For the proof of this theorem we use again the identification of $\mathcal{M}$ with $\mathrm{PSL}(2,\mathbb{R})$. The projective 
cocycle then becomes a matrix cocycle, for which we can apply the criterium of \cite{bochi-navas} (see also \cite{bochi-navas-2}). 
According to this, it suffices to establish that, for each element $f$ in the acting group, the Lyapunov exponent of this cocycle 
vanishes. In concrete terms, this amounts to saying that the parameters 
$\kappa_n =\kappa_n(\theta)$ and $\sigma_n = \sigma_n (\theta)$ involved in the expression of 
$P_{f^n,\theta}$, namely
$$e^{i P_{f^n,\theta} (\xi)} = e^{i \kappa_n} \cdot \frac{e^{i \xi} - \sigma_n}{1 - \bar{\sigma}_n e^{i \xi}},$$  
are such that the sequence 
\begin{equation}\label{drift}
\frac{\mathrm{dist}_{\mathbb{D}} (0,\sigma_n)}{n}
\end{equation}
converges to zero as $n \to \infty$ uniformily in $\theta$. 
Here, $\mathrm{dist}_{\mathbb{D}}$ stands for the hyperbolic distance, hence 
$$\mathrm{dist}_{\mathbb{D}} (0,\sigma_n) = \log \left( \frac{1 - r_n}{1 + r_n} \right),$$
where $r_n := \|\sigma_n||$.

To check the linear drift above, remind from (\ref{erre}) that 
\begin{equation}\label{erre-n}
r_n = \sqrt{\frac{\mathrm{D}_n^2 (\mathrm{D}_n - 1)^2 + \Delta_n^2}{\mathrm{D}_n^2 (\mathrm{D}_n + 1)^2 + \Delta_n^2}},
\end{equation}
where $\mathrm{D}_n = Df^n (\theta)$ and $\Delta_n = D^2 f^n (\theta)$. Now, since our group acts freely by $C^2$ diffeomorphisms 
and is finitely generated,  
its action is topologically conjugate to that of a group of rotations \cite{book}. Therefore, if $f \in \Gamma$ has rational 
rotation number, then it has finite order, hence the sequence (\ref{drift}) obviously converges to zero. If $f \in \Gamma$ 
has irrational rotation number, then it is well known that the growth of the derivatives along its iterates is uniformly 
subexponential: for every $\mu > 0$ there exists $C > 0$ such that, at every point $\theta$, one has
\begin{equation}\label{derivatives}
\frac{e^{- n \mu } }{C} \leq \mathrm{D}_n \leq C \, e^{n \mu }.
\end{equation}
Moreover, the cocycle rule for the affine derivative yields, for $M := \max | D^2 f / Df |$, 
$$\Delta_n 
= D f^n \cdot \frac{D^2 f^n}{D f^n} 
= D f^n \sum_{k=0}^{n-1} \frac{D^2 f}{Df} (f^k) \cdot (Df^k) 
\leq C e^{n \mu} \sum_{k=0}^{n-1} M C e^{k \mu} 
\leq C' e^{2n \mu},$$
where $C' = C' (C, M, \mu)$. Together with (\ref{erre-n}) and (\ref{derivatives}), this yields 
\begin{eqnarray*}
\frac{\mathrm{dist}_{\mathbb{D}} (0,\sigma_n)}{n} 
&=& 
\frac{1}{n} \log \left( \frac{1 + r_n}{1 - r_n} \right) \\ 
&=& 
\frac{1}{n} \log \left( \frac{(1 + r_n)^2}{1 - r_n^2} \right) \\
&\leq&
\frac{\log(4)}{n} + \frac{1}{n} \log \left( \frac{1}{1 - r_n^2} \right) \\
&=&
\frac{\log(4)}{n} + \frac{1}{n} \log \left( \frac{\mathrm{D}_n^2 (\mathrm{D}_n + 1)^2 + \Delta_n^2}{4 \mathrm{D}_n^3} \right) \\
&=&
\frac{\log(4)}{n} + \frac{1}{n} \log \left( \frac{\mathrm{D}_n^4 + 2\mathrm{D}_n^3 + \mathrm{D}_n^2 + \Delta_n^2}{4 \mathrm{D}_n^3} \right) \\
&\leq&
\frac{\log(4)}{n} + \frac{1}{n} \max \left\{ \log \left( \mathrm{D}_n \right), \log \left( 2 \right), 
\log \left( \frac{1}{\mathrm{D}_n} \right), \log \left( \frac{\Delta_n^2}{\mathrm{D}_n^3} \right) \right\} \\
&\leq&
\frac{\log(4)}{n} + \frac{1}{n} \max \left\{ \log \left( C e^{n \mu} \right), \log (2), \log \left( \frac{C^3 C'^2 e^{2 n \mu}}{e^{-3 n \mu}} \right) \right\}. 
\end{eqnarray*}
For large enough $n$, the last expression is smaller than or equal to \, $6 \log (\mu)$. \, 
Since this holds for every $\mu > 0$,  the sequence (\ref{drift}) uniformly converges to zero, which concludes the proof of the theorem.

\vspace{0.2cm}

\begin{obs}
{\em It is still unknown whether for each $C^2$ circle diffeomorphism $f$ with irrational rotation number there exists a sequence of $C^2$ diffeomorphisms 
$\varphi_n$ such that the conjugates $\varphi_n f \varphi_n^{-1}$ converge to a rotation as $n \to \infty$. (This question stands more generally for groups 
of $C^2$ diffeomorphisms acting freely; see \cite{Na:first}.) If this was the case, then a natural conjugating map $B_n$ for the projective cocycle $P_{f,\cdot}$ 
would be $B_n = P_{\varphi_n,\theta}$, because of the cocycle relation (compare (\ref{eq:conjj})). 

Actually, the main motivation of this work was to use the projective cocycle to deal with the conjugacy question above. However, our positive answer is only 
partial (at the level of cocycles), since there is no reason to expect that the conjugating cocycles $B_n$ we prove to exist arise from circle diffeomorphisms 
$\varphi_n$.}
\end{obs}

\vspace{0.2cm}

We next provide a partial converse to Theorem B.
Notice that the case of finite orbits 
that we avoid in the statement below essentially reduces to that of groups of diffeomorphisms of the interval, 
for which the behavior of the projective cocycle seems hard to understand (compare \cite{eynard-yo}).

\vspace{0.2cm}

\begin{prop} \label{p:converse}
Let $\Gamma$ be a finitely generated group of $C^2$ circle diffeomorphisms having no finite orbit. 
If the projective cocycle over its action is almost $\mathrm{SO} (2,\mathbb{R})$-reducible, then $\Gamma$ is Abelian, 
and its action is topologically conjugate to that of a group of rotations.
\end{prop}

\noindent{\bf Proof.} The converse of the criterium of \cite{bochi-navas,bochi-navas-2} also holds. (Actually, this is the 
easy-to-establish part of the criterium.) Therefore, to prove the proposition, it suffices to show that, if the group is not  topologically 
conjugate to a group of rotations, then it contains elements along whose iterates the drift (\ref{drift}) is positive. To do this, recall 
that Sacksteder's theorem and its generalizations imply that, if a finitely-generated 
group of $C^2$ circle diffeomorphisms is not topologically conjugate 
to a group of rotations, then it contains elements with hyperbolic fixed points  \cite{DKN-acta}. The proof will be hence finished by 
establishing that the drift (\ref{drift}) along the iterates of a diffeomorphism with hyperbolic fixed points is positive. 

Let $\theta_0$ be a point such that $f (\theta_0) = \theta_0$ and $Df (\theta_0) = e^{\mu} > 1$ for a $C^2$ circle diffeomorphism $f$. 
At this point we have $\mathrm{D}_n = e^{n \mu}$ for all $n$. Therefore,
\begin{eqnarray*}
\frac{\mathrm{dist}_{\mathbb{D}} (0,\sigma_n)}{n}  
=
\frac{1}{n} \log \left( \frac{1 + r_n}{1 - r_n} \right) 
=
\frac{1}{n} \log \left( \frac{(1 + r_n)^2}{1 - r_n^2} \right) 
&\geq &
\frac{1}{n} \log \left( \frac{1}{1 - r_n^2} \right) \\
&=&
\frac{1}{n} \log \left( \frac{\mathrm{D}_n^2 (\mathrm{D}_n + 1)^2 + \Delta_n^2}{4 \mathrm{D}_n^3} \right) \\
&\geq&
\frac{1}{n} \log \left( \frac{\mathrm{D}_n}{4} \right) \\
&=&
\frac{1}{n} [n \mu - \log (4)].
\end{eqnarray*}
This shows that the drift is larger than or equal to $\mu$, as desired.  $\hfill\square$


\section*{Appendix: The enlarged projective cocycle}

For a circle homeomorphism $f$ and three different angles $\theta_1,\theta_2,\theta_3$, we let $M_{f,(\theta_1,\theta_2,\theta_3)}$ 
be the unique M\"obius transformation that sends each $\theta_i$ into $f(\theta_i)$. It is straightforward to check that, in fractional notation,
\begin{equation}\label{e:oto}
M_{f,(\theta_1,\theta_2,\theta_3)} (z) = \frac{z (C - A \tau)  + (A c \tau - a C) }{z (1 - \tau)  + (c \tau - a)},
\end{equation}
where 
$$a := e^{i\theta_1}, \,\, b := e^{i \theta_2}, \,\, c := e^{i \theta_3}, \,\, 
A := e^{i f (\theta_1)}, \,\, B := e^{i f (\theta_2)}, \,\, C := e^{i f (\theta_3)}, \,\, 
\tau := \frac{b-a}{B-A} \cdot \frac{C-B}{c-b}.$$
It is better to write expression (\ref{e:oto}) in the form
\begin{equation}\label{e:normal-form}
M_{f,(\theta_1,\theta_2,\theta_3)}(z) 
= \frac{C - A \tau}{c \tau - a} \left( \frac{z + \frac{A c \tau - a C}{C - A \tau} }{1 + z \left( \frac{1 - \tau}{c\tau - a} \right) } \right).
\end{equation}
Indeed, as one can easily check, this corresponds to the canonical form of a M\"obius transformation
$$z \mapsto e^{i\kappa} \cdot \frac{z - \sigma}{1-\bar{\sigma} z},$$
with 
\begin{equation}\label{e:coeff}
e^{i \kappa} 
=  \frac{C - A \tau}{c \tau - a} , \qquad \sigma 
= -  \frac{A c \tau - a C}{C - A \tau} = - \overline{\left( \frac{1-\tau}{c\tau-a} \right)} = r e^{i\alpha}.
\end{equation}

By the definition, the following cocycle relation holds:
$$M_{f,(g(\theta_1),g(\theta_2),g(\theta_3))} \circ M_{g,(\theta_1,\theta_2,\theta_3)} = M_{fg,(\theta_1,\theta_2,\theta_3)}.$$
In particular, 
$$M_{f^{-1},(\theta_1,\theta_2,\theta_3)} = M_{f,(f^{-1}(\theta_1),f^{-1}(\theta_2),f^{-1}(\theta_3))}^{-1}.$$
We will call $M$ the {\em enlarged projective cocycle}. The justification of this terminology comes from the next proposition.

\vspace{0.1cm}

\begin{prop} \label{prop:gato}
If $f$ is a $C^2$ diffeomorphism, then $M_{f,\cdot}$ extends continuously to the whole product space $(\mathbb{S}^1)^3$. 
Moreover, the restriction of this extension to the diagonal coincides with the projective cocycle. 
\end{prop}

\noindent{\bf Proof.} We need to prove that the coefficients in expression (\ref{e:normal-form}), namely those in 
(\ref{e:coeff}), converge to those of (\ref{e:gen-normal}). We will perform the explicit computations in the case of triplets of 
different points converging to a triplet in the diagonal, which is the most interesting one. (The other cases are left to the reader.) 
We hence let $\theta_1, \theta_2, \theta_3$ be three distinct cycically ordered points on the circle converging to the same $\theta$.
For simplicity, we let $\varepsilon := (\theta_3 - \theta_1)/2$. Accordingly, and with some abuse of notation, parameters arising 
in (\ref{e:normal-form}) and (\ref{e:coeff}) will be denoted with a subindex $\varepsilon$.

Let us write $\tau_{\varepsilon} = \rho_{\varepsilon} e^{i \gamma_\varepsilon}$, where $\rho_\varepsilon := \| \tau_\varepsilon \|$. Since 
$$\rho_\varepsilon e^{-\gamma_\varepsilon} 
= \overline{\tau}_\varepsilon 
= {\frac{\bar{b}_\varepsilon-\bar{a}_\varepsilon}{\bar{B}_\varepsilon-\bar{A}_\varepsilon} \cdot \frac{\bar{C}_\varepsilon-\bar{B}_\varepsilon}{\bar{c}_\varepsilon-\bar{b}_\varepsilon}} 
= {\frac{1 / b_\varepsilon - 1 / a_\varepsilon}{1 / B_\varepsilon - 1 / A_\varepsilon} \cdot \frac{1 / C_\varepsilon - 1 / B_\varepsilon}{1 / c_\varepsilon - 1 / b_\varepsilon}} 
= \frac{A_\varepsilon c_\varepsilon}{a_\varepsilon C_\varepsilon} \tau_\varepsilon 
= \frac{A_\varepsilon c_\varepsilon}{a_\varepsilon C_\varepsilon} \rho_\varepsilon e^{\gamma_\varepsilon},$$
we have
$$\gamma_\varepsilon = \frac{f (\theta_3) - f (\theta_1)}{2} -  \frac{\theta_3 - \theta_1}{2}.$$ 
Using this, one computes 
$$
\| \sigma_\varepsilon \|^2 = \left\| \frac{1-\tau}{c\tau-a} \right\|^2  
=
\frac{1 - 2 \, \rho_\varepsilon \cos(  \gamma_{\varepsilon} ) + \rho_\varepsilon^2}
{1 - 2 \, \rho_\varepsilon \cos(  \gamma_{\varepsilon} + \theta_3 - \theta_1 ) + \rho_\varepsilon^2} 
=
\frac{ 2  \rho_\varepsilon ( 1 - \cos(  \gamma_{\varepsilon} ))  + (1 -  \rho_\varepsilon)^2}
{2  \rho_\varepsilon ( 1 - \cos(  \gamma_{\varepsilon} + \theta_3 - \theta_1 )) + ( 1 - \rho_\varepsilon)^2} .
$$
We claim that the following convergences hold as $\varepsilon \to 0$:
\begin{small}
\begin{equation}\label{e:tres}
\frac{1-\cos (\gamma_{\varepsilon})}{\varepsilon^2} \to \frac{( Df(\theta) - 1)^2}{2}, 
\quad
\frac{1-\cos (\gamma_{\varepsilon} + \theta_3 - \theta_1 )}{\varepsilon^2} \to \frac{( Df(\theta) + 1)^2}{2}, 
\quad 
\frac{1 - \rho_{\varepsilon}}{\varepsilon} \to - \frac{D^2 f(\theta)}{Df (\theta)}.
\end{equation}
\end{small}The first two convergences in (\ref{e:tres}) follow analogously just by noticing that 
$$\gamma_{\varepsilon} + \theta_3 - \theta_1 = \frac{f(\theta_3)-f(\theta_1)}{2} + \frac{\theta_3-\theta_1}{2}.$$
We thus consider only the first one, which is shown as follows:
\begin{eqnarray*}
\frac{1\!-\!\cos (\gamma_{\varepsilon})}{\varepsilon^2} 
&=& \frac{\frac{\gamma_{\varepsilon}^2}{2} + o(\varepsilon^2)}{ \varepsilon^2} \\
&=& \frac{( \frac{f(\theta_3) - f(\theta_1))}{2} - \frac{\theta_3 - \theta_1}{2} )^2 + o(\varepsilon^2)}{2 \varepsilon^2} \\
&=& \frac{(\frac{Df(\theta_1) (\theta_3-\theta_1)}{2} - \frac{\theta_3 - \theta_1}{2})^2 + o(\varepsilon^2)}{2 \varepsilon^2} \\
&=& \frac{(Df(\theta_1 ) - 1)^2 + o(\varepsilon^2)}{2\varepsilon^2} \quad \longrightarrow \quad  \frac{( Df(\theta) - 1 )^2}{2}.
\end{eqnarray*}
The third expression in (\ref{e:tres}) is more interesting: 
\begin{eqnarray*}
\frac{1 - \rho_{\varepsilon}}{\varepsilon} 
&=& \frac{1 - \rho_{\varepsilon}^2}{(1 + \rho_{\varepsilon}) \, \varepsilon} 
\,\, = \, \, \frac{1 - \rho_{\varepsilon}^2}{( 2 + o(1) ) \, \varepsilon} 
\,\, =  \,\, \frac{ 1 - \frac{ \| C_{\varepsilon} - B_{\varepsilon} \|^2 }{ \| c_{\varepsilon} - b_{\varepsilon} \|^2 } \Big/ 
 \frac{ \| B_{\varepsilon} - A_{\varepsilon} \|^2 }{ \| b_{\varepsilon} - a_{\varepsilon} \|^2}  }
{2 \, \varepsilon + o (\varepsilon)} 
\,\, = \,\, \frac{  \frac{ \| B_{\varepsilon} - A_{\varepsilon} \|^2 }{ \| b_{\varepsilon} - a_{\varepsilon} \|^2} - 
  \frac{ \| C_{\varepsilon} - B_{\varepsilon} \|^2 }{ \| c_{\varepsilon} - b_{\varepsilon} \|^2} }
{2 \,  \frac{ \| B_{\varepsilon} - A_{\varepsilon} \|^2 }{ \| b_{\varepsilon} - a_{\varepsilon} \|^2}  \, \varepsilon + o (\varepsilon) }  \\
&=&
\frac{ \frac{1 - \cos (f(\theta_2)-f(\theta_1))}{1 - \cos (\theta_2 - \theta_1)} - 
 \frac{1 - \cos (f(\theta_3)-f(\theta_2))}{1 - \cos (\theta_3 - \theta_2)} }
{2 \,  \frac{ \| B_{\varepsilon} - A_{\varepsilon} \|^2 }{ \| b_{\varepsilon} - a_{\varepsilon} \|^2}  \, \varepsilon + o (\varepsilon)}  
\,\, = \,\, 
\frac{ \frac{ ( f(\theta_2)-f(\theta_1))^2 }{ (\theta_2 - \theta_1)^2} - 
 \frac{( f(\theta_3)-f(\theta_2) )^2}{ (\theta_3 - \theta_2 )^2} + o (\varepsilon^2)}
{2 \,  Df(\theta)^2 \, \varepsilon + o (\varepsilon)} \\
&=&
\frac{ \left( Df(\theta_2) - \frac{D^2f (\theta_2)}{2} (\theta_2 - \theta_1) \right)^2 - 
 \left( Df(\theta_2) + \frac{D^2f (\theta_2)}{2} (\theta_3 - \theta_2) \right)^2 + o (\varepsilon^2)}
{2 \,  Df(\theta)^2 \, \varepsilon + o (\varepsilon)} \\ 
&=&
- \frac{ Df(\theta_2) \, D^2 f (\theta_2) \, ( (\theta_2-\theta_1) +(\theta_3-\theta_2)) + o (\varepsilon^2)}
{2 \,  Df(\theta)^2 \, \varepsilon + o (\varepsilon)}  \\
&=& - \frac{ Df(\theta_2) \, D^2 f (\theta_2) \,  2 \varepsilon + o (\varepsilon^2)}
{2 \,  Df(\theta)^2 \, \varepsilon + o (\varepsilon)} \quad
\longrightarrow  \quad -\frac{D^2f(\theta)}{Df (\theta)}.
\end{eqnarray*}
Putted together, the convergences (\ref{e:tres}) yield 
\begin{equation}\label{e:erre-e}
\| \sigma_{\varepsilon} \|^2 = r_\varepsilon^2 \longrightarrow  
\frac{ ( Df(\theta)^2 \, ( Df(\theta) - 1 )^2 + ( D^2 f (\theta) )^2}{ ( Df(\theta)^2 \, ( Df(\theta) + 1 )^2  + ( D^2 f (\theta) )^2},
\end{equation}
which is compatible with (\ref{r_sub_ene}). 

\vspace{0.1cm}

Let us now look at the angle $\alpha_\varepsilon$ of $\sigma_{\varepsilon} = r_\varepsilon e^{i \alpha_\varepsilon}$. 
First notice that, from (\ref{e:coeff}), 
$$r_{\varepsilon} e^{-i \alpha_{\varepsilon}} 
= \bar{\sigma}_{\varepsilon} = \frac{\tau_{\varepsilon} - 1}{c_{\varepsilon} \tau - a_{\varepsilon}} 
= \frac{ \tau_{\varepsilon} - 1 }{ a_{\varepsilon} ( \frac{c_{\varepsilon} \tau_{\varepsilon} }{a_{\varepsilon}} - 1 ) },$$
which may be rewritten as
$$r_{\varepsilon} e^{i (\theta_1 - \alpha_{\varepsilon})} 
= \frac{ \tau_{\varepsilon} - 1 }{  \frac{c_{\varepsilon} \tau}{a_{\varepsilon}}-1 } 
= \frac{ (\rho_{\varepsilon} \cos (\gamma_{\varepsilon}) - 1) 
+ i \rho_{\varepsilon} \sin (\gamma_{\varepsilon})}{( \rho_{\varepsilon} \cos (\gamma_{\varepsilon} + \theta_3 - \theta_1) - 1 ) 
+ i \rho_{\varepsilon} \sin (\gamma_{\varepsilon} + \theta_3 - \theta_1) }.$$
A straightforward computation then yields 
\begin{equation}\label{e:tan}
\tan(\theta_1 - \alpha_{\varepsilon}) 
= \frac{\rho_{\varepsilon} \sin (\gamma_{\varepsilon}) ( \rho_{\varepsilon} \cos (\gamma_{\varepsilon} + \theta_3 - \theta_1) - 1)  
- ( \rho_{\varepsilon} \cos (\gamma_{\varepsilon}) - 1 ) \rho_{\varepsilon} \sin (\gamma_{\varepsilon} + \theta_3 - \theta_1)}
{ ( \rho_{\varepsilon} \cos (\gamma_{\varepsilon}) - 1)( \rho_{\varepsilon} \cos (\gamma_{\varepsilon} + \theta_3 - \theta_1) - 1) 
+ \rho_{\varepsilon} \sin (\gamma_{\varepsilon}) \rho_{\varepsilon} \sin (\gamma_{\varepsilon} + \theta_3 - \theta_1)}.
\end{equation}
Using (\ref{e:tres}) together with 
$$\sin (\gamma_{\varepsilon}) = \varepsilon \, ( Df(\theta) - 1) + o (\varepsilon^2), 
\qquad \quad \sin (\gamma_{\varepsilon} + \theta_3 - \theta_1) = \varepsilon \, (Df(\theta) + 1) + o (\varepsilon^2),$$
the right-side expression in (\ref{e:tan}) is easily proven to be of the order
$$\frac
{\varepsilon^2 \left(  \frac{D^2 f (\theta)(Df (\theta) - 1 )}{Df(\theta)} - \frac{D^2 f (\theta)(Df (\theta) + 1)}{Df(\theta)} \right) + o (\varepsilon^2)}
{\varepsilon^2 \left( (\frac{D^2f(\theta)}{Df(\theta)})^2 + (Df(\theta)^2 - 1) \right) + o (\varepsilon^2)}.$$
Passing to the limit as $\varepsilon \to 0$, we thus deduce that 
$$\tan (\theta - \alpha) = - \frac{2 \, Df (\theta) \, D^2 f (\theta)}{Df(\theta)^2 ( Df (\theta)^2 - 1) + (D^2f(\theta))^2}$$
whenever the denominator above is nonzero, which is compatible with (\ref{e:angle-sigma}). This shows that, in this case, $\alpha_{\varepsilon} \to \alpha$, 
hence, $\sigma_\varepsilon \to \sigma$. In case of zero denominator (that is, when $Df(\theta) = 1$ and $D^2 f (\theta) = 0$), relation (\ref{e:erre-e}) 
shows that $\sigma_{\varepsilon} \to 0$, which still corresponds to  $\sigma_\varepsilon \to \sigma$.

\vspace{0.1cm}

We have thus proved that $\sigma_\varepsilon$ converges to the parameter $\sigma$ from (\ref{e:gen-normal}). 
Using this fact and
$$e^{i\kappa_{\varepsilon}} \cdot \frac{e^{i \theta_2} - \sigma_{\varepsilon}}{1 - \bar{\sigma}_{\varepsilon} e^{i \theta_2}},
= e^{iM_{f,(\theta_1,\theta_2,\theta_3)}(\theta_2)} 
= e^{i P(\theta_2)} 
= e^{i\kappa} \cdot \frac{e^{i \theta_2} - \sigma}{1 - \bar{\sigma} e^{i \theta_2}},$$
we obtain that $e^{i\kappa_{\varepsilon}} \to e^{i \kappa}$. Therefore, $\kappa_{\varepsilon} \to \kappa$, thus closing the proof.
$\hfill\square$


\vspace{0.25cm}

We next state an extension of Theorem A to the cocycle $M$.

\vspace{0.5cm}

\noindent{\bf Theorem A'}
{\em The enlarged projective cocycle associated to a group of $C^2$ circle diffeomorphisms is 
$\mathrm{SO} (2,\mathbb{R})$-reducible if and only if the action is $C^2$ conjugate to an action by rotations.}


\vspace{0.5cm}

Indeed, if $\varphi$ is a $C^2$ diffeomorphism that conjugates the action to that of a group of rotations, then 
$(\theta_1,\theta_2,\theta_3) \to M_{\varphi, (\theta_1,\theta_2,\theta_3)}$ performs the reduction of the cocycle 
$M$ into a cocycle of rotations. Conversely, if $M$ is $\mathrm{SO} (2,\mathbb{R})$-reducible, then its restriction 
to the diagonal is also $\mathrm{SO} (2,\mathbb{R})$-reducible. By Proposition \ref{prop:gato} and Theorem A, 
the group is $C^2$-conjugate to a group of rotations.

\vspace{0.25cm}

We close with an extension of Theorem B.

\vspace{0.4cm}

\noindent{\bf Theorem B'} {\em If $\Gamma$ is a finitely generated Abelian group of $C^2$ circle diffeomorphisms acting 
freely, then the enlarged projective cocycle above its action is almost $\mathrm{SO}(2,\mathbb{R})$-reducible.}

\vspace{0.4cm}

Notice that a partial converse holds for this result: if $\Gamma$ is a finitely generated group of $C^2$ circle diffeomorphisms having 
no finite orbit for which the enlarged projective cocycle is almost $\mathrm{SO} (2,\mathbb{R})$-reducible, then $\Gamma$ is Abelian. Indeed, 
it suffices to apply Proposition \ref{p:converse} to the restriction of $M$ to the diagonal to obtain this conclusion.

Unfortunately, Theorem B' does not follow as a direct consequence of Theorem B, though the schema of proof is similar. Namely, 
following again \cite{bochi-navas,bochi-navas-2},
we need to check that, for each element $f$ in the acting group, the parameters 
$\kappa_n =\kappa_n(\theta_1,\theta_2,\theta_3)$ and $\sigma_n = \sigma_n (\theta_1,\theta_2,\theta_3)$ involved in the expression of 
$M = M_{f^n,(\theta_1,\theta_2,\theta_3)}$, namely
$$M (z) = e^{i \kappa_n} \cdot \frac{z - \sigma_n}{1 - \bar{\sigma_n} z},$$  
are such that the sequence $\frac{1}{n} \log (\frac{1}{1 - r_n^2})$ converges to zero as $n \to \infty$ uniformly in 
$(\theta_1,\theta_2,\theta_3)$, where $r_n := \| \sigma_n \|$. We sketch the computations below in the most relevant case, 
namely when $\theta_1,\theta_2,\theta_3$ are all distinct and cyclically ordered on the circle (for reverse ordered triplets, 
just interchange $\theta_1$ and $\theta_2$). The reader is invited to check that the estimates we obtain are uniform in 
$(\theta_1,\theta_2,\theta_3)$, which is necessary to complete the proof for the general case.
 
As in the proof of Proposition \ref{prop:gato}, one easily computes
\begin{equation}\label{e:expression}
\frac{1}{1 - r_n^2} =
\frac{(1 - \rho_n)^2 + 2 \, \rho_n (1 - \cos(  \frac{\theta_3 - \theta_1}{2} + \frac{f^n(\theta_3) - f^n(\theta_1)}{2} ))}
{4 \, \rho_n \sin( \frac{\theta_3 - \theta_1}{2} ) \sin(\frac{f^n(\theta_3) - f^n(\theta_1)}{2} )},
\end{equation}
where $\rho_n = \|\tau_n\|$. 
Since the group we are dealing with is topologically conjugate to a group of rotations, its action on the circle is equicontinuous. 
We can hence fix small-enough positive constants $\delta, \bar{\delta}$ such that $| g(\bar{\theta}) - g(\theta) | \geq \delta$ for 
all $g$ in the acting group whenever $| \bar{\theta} - \theta| \geq \bar{\delta}.$ In order to estimate expression (\ref{e:expression}), 
there are hence two regimes to consider.

\vspace{0.2cm}

\noindent{$\bullet$ \bf Assume $|\theta_3 - \theta_1| \geq \bar{\delta}$.}

\vspace{0.2cm}

Letting $s = \sin (\delta/2)$, we get
$$\frac{1}{1 - r_n^2} \leq \frac{1 + \rho_n^2}{4 \, \rho_n \, s^2}.$$
By the mean value theorem,
$$\rho_n = 
\left\|
\frac{
\frac{e^{if^n(\theta_3)}-e^{if^n(\theta_2)}}{e^{i\theta_3}-e^{i\theta_2}}  
}
{
\frac{e^{if^n(\theta_2)}-e^{if^n(\theta_1)}}{e^{i \theta_2}-e^{i \theta_1}} } 
\right\|$$
behaves as \, $D f^n (\xi_2) / D f^n (\xi_1) $ \, for certain points $\xi_1,\xi_2$.  
Since the growth (and decay) of derivatives along iterates of $f$ is uniformly subexponential, this allows easily establishing that 
$1/(1-r_n^2)$ is uniformly subexponential in this case. 

\vspace{0.2cm}

\noindent{$\bullet$ \bf Assume $|\theta_3 - \theta_1| \leq \bar{\delta}$.}

\vspace{0.2cm}

In this case, one also has $|g(\theta_3) - g(\theta_1)| \leq \hat{\delta}$ for a certain small constant $\hat{\delta}$ and all $g$ in the acting group, 
due to equicontinuity. We will decompose expression (\ref{e:expression}) in two terms.

\vspace{0.1cm}

The first term is
\begin{equation}\label{e:first}
\frac{(1 - \rho_n)^2 }
{2 \, \rho_n \sin( \frac{\theta_3 - \theta_1}{2} ) \sin(\frac{f^n(\theta_3) - f^n(\theta_1)}{2} )}.
\end{equation}
Since both $|\theta_3-\theta_1|$ and $|f^n(\theta_3)-f^n(\theta_1)|$ are small, $\sin (\frac{\theta_3-\theta_1}{2})$ is equivalent to 
$|\theta_3-\theta_1|$, and $\sin (\frac{f^n(\theta_3).- f^n(\theta_1}{2})$ is equivalent to 
$|f^n(\theta_3) - f^n(\theta_1)| = |Df^n (\xi) (\theta_3 - \theta_1)|$ for a certain point $\xi$. 
Since $\rho_n$ and $Df^n$ have subexponential behavior, 
it remains to compare $|1-\rho_n|$ with $|\theta_3 - \theta_1|$. Now, this corresponds to the last expression in (\ref{e:tres}). 
By performing similar estimates but taking care now of the order of the coefficients involved in the $o (\varepsilon^2)$ terms, 
one easily deduces that the quotient $|1-\rho_n| / |\theta_3-\theta_1|$ 
is controlled by the first and second derivatives of $f$, which behave subexponentially.

\vspace{0.1cm}

The second term to analyze is 
\begin{equation}\label{e:second}
\frac{2 \, \rho_n (1 - \cos(  \frac{\theta_3 - \theta_1}{2} + \frac{f^n(\theta_3) - f^n(\theta_1)}{2} ))}
{2 \, \rho_n \sin( \frac{\theta_3 - \theta_1}{2} ) \sin(\frac{f^n(\theta_3) - f^n(\theta_1)}{2} )} 
=
\frac{1 - \cos(  \frac{\theta_3 - \theta_1}{2} + \frac{f^n(\theta_3) - f^n(\theta_1)}{2} )}
{ \sin( \frac{\theta_3 - \theta_1}{2} ) \sin(\frac{f^n(\theta_3) - f^n(\theta_1)}{2} )}. 
\end{equation}
This is of the order 
$$\frac{\frac{1}{2} (\frac{\theta_3 - \theta_1}{2} + \frac{f^n(\theta_3) - f^n(\theta_1)}{2})^2}{\frac{|\theta_3 - \theta_1|}{2} \cdot \frac{|f^n(\theta_3) - f^n (\theta_1)|}{2}}  
= \frac{|\theta_3 - \theta_1|^2 (1+Df^n (\hat{\xi}))^2}{2 \, |\theta_3 - \theta_1|^2 \, Df^n (\xi)}
= \frac{ (1+Df^n (\hat{\xi}))^2}{2 \, Df^n (\xi)}.$$
Again, the subexponential behavior of $Df^n$ allows us to conclude.

\vspace{0.5cm}

\noindent{\bf Acknowledgments.} Andr\'es Navas was funded by the projects FONDECYT 1200114 (in Chile) as well as
 FORDECYT 265667 and the PREI of the DGAPA at UNAM (in M\'exico). Mario Ponce was funded by the projects 
 FONDECYT 1180922 as well as ANILLO ACT172001 CONICYT.


\begin{small}

\vspace{0.2cm}

\noindent Andr\'es Navas\\
\noindent Dpto de Matem\'atica y C.C., Universidad de Santiago de Chile\\
\noindent Alameda 3363, Estaci\'on Central, Santiago, Chile\\
\noindent and\\
\noindent Unidad Cuernavaca Instituto de Matem\'aticas\\
\noindent Universidad Nacional Aut\'onoma de M\'exico, Campus Morelos\\
\noindent E-mail: andres.navas@usach.cl\\
\vspace{0.2 cm}

\noindent Mario Ponce\\
\noindent Facultad de Matem\'aticas, Pontificia Universidad Cat\'olica de Chile\\
\noindent Av. Vicu\~na Mackenna 4860, Macul, Chile\\
\noindent E-mail: mponcea@mat.uc.cl\\

\end{small}


\begin{thebibliography}{Dillo 83}


\bibitem{bochi-navas} {\sc J. Bochi \& A. Navas.} A geometric path from zero Lyapunov 
exponents to rotation cocycles. {\em Erg. Theory and Dyn. Systems} {\bf 35} (2015), 374-402.

\bibitem{bochi-navas-2} {\sc J. Bochi \& A. Navas.} Almost reduction and perturbation of matrix cocycles. 
{\em Annales de l'Institut Herni Poincar\'e Analyse non Lin\'eaire} {\bf 6} (2014), 1101-1107.

\bibitem{CNP}{\sc D. Coronel, A. Navas \& M. Ponce.} On bounded cocycles of isometries over a minimal dynamics. 
{\em J. Mod. Dyn.} {\bf 7} (2013),  45-74. 

\bibitem{DKN-acta} {\sc B. Deroin, V. Kleptsyn \& A. Navas.} 
Sur la dynamique unidimensionnelle en r\'egularit\'e interm\'ediaire. {\em Acta Math.}  {\bf 199} (2007), 199-262.

\bibitem{eynard-yo} {\sc H. Eynard \& A. Navas.} 
Mather invariant, distortion, and conjugates for diffeomorphisms of the interval. 
Preprint (2019); arXiv:1912.09305.

\bibitem{GALA}{\sc F. Gardiner \& N. Lakic.} {\em Quasiconformal Teichm\"uller Theory.} Mathematical 
Surveys and Monographs Vol. {\bf 76}. American Mathematical Society, Providence, RI, (2000).

\bibitem{herman} {\sc M. Herman.} Sur la conjugaison diff\'erentiable des diff\'eomorphismes du cercle 
\`a des rotations.  {\em Inst. Hautes \'Etudes Sci. Publ. Math.} {\bf 49} (1979), 5-233. 



\bibitem{Na:first} {\sc A. Navas.} On conjugates and the asymptotic distortion of 1-dimensional $C^{1+bv}$ 
diffeomorphisms. Preprint (2018); arXiv:1811.06077.
 
\bibitem{compositio} {\sc A. Navas.} Sur les rapprochements par conjugaison en dimension 1 et classe $C^1$.  
{\em Compos. Math.} {\bf 150} (2014), 1183-1195. 

\bibitem{book} {\sc A. Navas.} {\em Groups of circle diffeomorphisms.}
Chicago Lect. in Mathematics (2011).

\bibitem{3-remarks} {\sc A. Navas.} Three remarks on one-dimensional bi-Lipschitz conjugacies. 
Unpublished Note (2006); arXiv:0705.0034.


\bibitem{SING}{\sc D. Singer.} Diffeomorphisms of the circle and hyperbolic curvature. 
{\em Conform. Geom. Dyn.} {\bf 5} (2001), 1-5.	

\end{thebibliography}
\end{document}